\documentclass[11pt]{article}
\usepackage[inactive]{srcltx}
\usepackage{epsfig,amssymb}

\begin{document}

\title{
Note on generic singularities of planar flat 3-webs}

\author{{\Large Sergey I. Agafonov}\\
Departmento de Matem\'atica,\\
UNESP-Universidade Estadual Paulista,\\ S\~ao Jos\'e do Rio Preto, Brazil\\
e-mail: {\tt agafonov@ibilce.unesp.br} }
\date{}
\maketitle
\unitlength=1mm

\newtheorem{theorem}{Theorem}
\newtheorem{proposition}{Proposition}
\newtheorem{lemma}{Lemma}
\newtheorem{corollary}{Corollary}
\newtheorem{definition}{Definition}
\newtheorem{example}{Example}

\pagestyle{plain}

\begin{abstract}
\noindent We propose a definition of genericity for singular flat planar 3-webs formed by integral curves of implicit ODEs and give a classification of generic singularities of such webs. 
(secondary).
\end{abstract}

\section{Introduction}

 A planar $d$-web is formed by $d$ foliations in the plane. At each point of the plane, we have  $d$ leaves passing through the point, one from each foliation of the web.  A point is called {\bf regular} if any two of these $d$ leaves are transverse.
Consider the pseudogroup of local diffeomorphisms of the plane and the corresponding equivalence relation on the set of $d$-web germs. Any two planar 2-web germs are equivalent whenever the base points are regular.  This is not true for 3-web germs  (see \cite{BB}). There is a local invariant, which has, in fact, a topological nature. The differential-geometric counterpart of this invariant is the so-called Blaschke curvature.
\begin{definition}
3-web is flat (or hexagonal) if its germ at any regular point is equivalent to the web formed by three families of parallel
lines.
\end{definition}
A 3-web is flat if and only if the Blaschke curvature vanishes identically (see, for instance, \cite{BB}).
This curvature  is a scalar 2-form, therefore any general classification of 3-webs with respect to local diffeomorphisms
will necessarily have functional moduli. Namely, such a classification will inevitably involve arbitrary functions of two variables.

By the above definition, any two flat 3-web germs are equivalent provided that the base points are regular.
Hence the "personality" of a flat 3-web is encoded in its singularities.

Web structure is ubiquitous in mathematics and its applications, the Blaschke curvature often being the obstacle to obtaining a "reasonable" classification.
Therefore flat 3-webs play a distinguished role.

For example,
  hexagonal 3-webs have a 3-dimensional
local symmetry algebra at regular points, while a generic 3-web does not admit any infinitesimal symmetry (see  \cite{Cg}). (By infinitesimal symmetry of the web we understand a vector field whose local flow preserves the web.)

It was also observed that characteristic 3-webs of  integrable equations, playing important role in mathematical and theoretical physics, are flat (see, for instance, \cite{Fc,Fi,MFa,Al}).

The singularities of planar 2-webs, defined by solutions of implicit ODEs, quadratic in the derivative, are well understood (see \cite{Te,Dn,Ag}), whereas  singular implicit ODEs, polynomial in derivatives of degree $d\ge 3$ bring more difficulties, the main obstacle being the nontrivial web structure on its solutions (see \cite{Ds,Dn,Hs,Ns,Nw}).

Let us review some known results on classification of flat 3-webs.
One can show that a general classification of singular flat 3-webs will also have functional moduli, now the "moduli" being arbitrary functions of one variable (see, for example, the discussion in \cite{Ai}).
To get a reasonable classification one has to restrict the class of admissible singularities by imposing some meaningful conditions.

 A sufficiently general class of singular 3-web germs can be described by binary forms:
$$
K_3(x,y)dy^3+K_2(x,y)dy^2dx+K_1(x,y)dydx^2+K_0(x,y)dx^3=0.
$$
Dividing the above form by $dx^3$  one gets an implicit ODE, cubic in $p=\frac{dy}{dx}$:
\begin{equation}\label{cubic}
K_3(x,y)p^3+K_2(x,y)p^2+K_1(x,y)p+K_0(x,y)=0.
\end{equation}
Equation (\ref{cubic}) defines a (possibly singular) surface $M$ in 3-dimesional contact space $\mathbb R^2 \times  \mathbb P^1(\mathbb R)$.
It is immediate that if a point $(x,y)$ is not regular then the projection $\pi : M\to
\mathbb R^2,\  (x,y,p) \mapsto (x,y)$ is not a local diffeomorphism at least at one point $m$
 in the fiber $\pi^{-1}(x,y)$.

 Define the {\it criminant} $C$ of implicit ODE (\ref{cubic}) as the set of points on $M$, where the projection $\pi$ fails to be a local diffeomorphism, and the
{\it discriminant}  $\Delta$ of equation (\ref{cubic}) (or the {\it apparent
contour} of the surface $M$) as the image of the criminant under the projection: $\Delta=\pi(C)$.

Under some natural conditions of regularity for the surface $M$ and for the criminant $C$, implying the smoothness of $M$ and $C$, (see equation (\ref{regularity}) in the next section),  the following normal forms were obtained in \cite{Ai}:
\begin{equation}\label{theorem_callsification_regular}
1)\ p^3+px-y=0,\ \ \ \ \ \ \  2)\ p^3+2xp+y=0,\ \ \ \ \ \ \ 3)\ p^2=x,\ \ \ \ \ \ \ 4)\ p^2=y.
\end{equation}

\noindent {\bf Remark.} For the quadratic ODEs, the third root is
$\infty$.
The cubic normal forms were conjectured in \cite{Nw}.
\medskip

Another natural restriction on the singularity  is the existence of at least one infinitesimal symmetry at the singular point.
It turns out that at a singular point generators of a symmetry algebra became infinite or multi-valued, thus this restriction is not trivial (see \cite{As}).
We call a web homogeneous if it is
invariant with respect to the flow of  Euler vector field
\begin{equation}\label{euler}
E=c_1 x\frac{\partial}{\partial x}+c_2 y\frac{\partial}{\partial
y}, \ \ \ \ \ \ c_i={\rm const.}
\end{equation}
Classification of singular homogeneous 3-webs were obtained in \cite{Ac}.
To give the reader some idea about the functions emerging in this classification, we present here some normal forms:
$$
\begin{array}{l}
p^3+y^2p=\frac{2}{\sqrt{27}}y^3\tan(2\sqrt{3}x),\\
p^3+4x(y-\frac{4}{9}x^3)p+y^2+\frac{64}{81}x^6-\frac{32}{9}yx^3=0,\\
p^3+y^{3}p=y^{\frac{9}{2}}U\left(xy^{\frac{1}{2}}\right),
\end{array}
$$
where $U$ is expressed in terms of  the Legendre functions $P^{\mu}_{\nu}(z),Q^{\mu}_{\nu}(z)$.
Note that all the forms (\ref{theorem_callsification_regular}) also are homogeneous.

In this note we address the problem of genericity for singularities of flat 3-webs. Our considerations are local, all the objects involved are smooth or analytic.

\section{Genericity via transversality}\label{section_genericity}

The 3-web, formed by integral curves of  equation (\ref{cubic}), is
flat if and only if $K=(K_3,K_2,K_1,K_0)$ satisfy a certain (rather involved) nonlinear
PDE of second order.
Following Blaschke \cite{Be}, let us normalize one-forms,  vanishing on the
web leaves, to satisfy the condition $
\sigma_1+\sigma_2 + \sigma _3=0.$ For example, one may choose
$$
\sigma_1=(p_2-p_3)(dy-p_1dx),\ \ \sigma_2=(p_3-p_1)(dy-p_2dx),\ \
\sigma_3=(p_1-p_2)(dy-p_3dx),\ \
$$
where $p_1,p_2,p_3$ are the roots of (\ref{cubic}) at a non-singular point
$(x,y)$.  The above normalization defines the so-called Chern
connection  form $\gamma$ by  $d\sigma_i=\gamma \wedge \sigma_i$. Its derivative $d\gamma $ is the Blaschke curvature. Thus the hexagonality of the web amounts to the equation $d\gamma =0$, which involves the second order derivatives of $p_i$. Manipulating elementary symmetric functions of $p_i$, one rewrites this equation in terms of $K$ and its derivatives up to the second order.

The study of generic properties of solutions to PDEs was mainly concentrated on equations of order one. These properties are relatively well understood in terms of Legendrian and Lagrangian singularities (see \cite{AGV,GZ}).

Our approach is based on the Thom ideas on genericity via transversality and is close to those of \cite{DDs,Dg}, though results obtained in these papers are not applicable here, as the PDE $d\gamma =0$ for $K$, being linear in 2nd derivatives, has vanishing coefficients of these derivatives. Thus the Cauchy-Kowalevskaya Theorem does not work and no existence theorem can help  us to analyze the space of solutions.

Consider the set $\mathcal{H}$ of map germs  $K:(\mathbb R^2,0) \to \mathbb R^4$ which solve the equation $d\gamma =0$. We adopt here a naive and old-fashioned point of view on perturbations, namely, we call $\widetilde{K}\in \mathcal{H}$ a perturbation of $K\in \mathcal{H}$ if, at the base point, the values of $\widetilde{K}$ and its first derivatives are close to the corresponding values of $K$. (Note that our knowledge of the subspace of $\mathcal{H}$, describing singular hexagonal 3-webs, is rather limited. For instance, it is not clear if higher derivatives of a perturbation $\widetilde{K}$ are close to those of $K$.)

For each  $(k_3,k_2,k_1,k_0)\in \mathbb K^4 \setminus \{(0,0,0,0)\}$  there are solutions $K\in \mathcal{H}$ such that $K(0,0)=(k_3,k_2,k_1,k_0)$. In fact, a cubic form $g(p)=k_3p^3+k_2p^2+k_1p+k_0$ can be brought by a M\"obius transform (i.e. by a transform that is fractional linear in p) to one of the following forms: 1) $p^3$ (triple root), 2) $p^2$ (double root), 3) $p(p-1)(p+1)$ (distinct roots).  Thus, applying a suitable local diffeomorphism, we get a desired solution from the normal form  2) of classification list (\ref{theorem_callsification_regular}), if $g(p)$ has a triple root;  or from the normal form 3), if $g(p)$ has a double root; or from the non-singular flat 3-web. Therefore for any $K\in \mathcal{H}$ there is a small perturbation $\widetilde{K}$ such that $\widetilde{K}(0,0)\ne (0,0,0,0)$.
This means that all web directions are well defined at each point.
Choosing local coordinates so that none of the web leaves at the origin is tangent to the $y$-axis, one ensures $K_3(0,0)\ne 0$. Dividing  (\ref{cubic}) by $K_3$, one reduces our equation to a monic one:
\begin{equation}\label{general_cub}
p^3+a(x,y)p^2+b(x,y)p+c(x,y)=0.
\end{equation}
Finally, "killing" the coefficient by $p^2$
  by a differential analog of Tschirnhausen transformation
\begin{equation}\label{chirn}
y=f(\tilde{x},\tilde{y}),\ x=\tilde{x},\ \ {\rm with} \ \
3f_{\tilde{x}}+a(\tilde{x},f)=0
\end{equation}
 (see \cite{Ai} for more detail), one arrives at the normal form
\begin{equation}\label{cub}
F(x,y,p)=p^3+A(x,y)p+B(x,y)=0.
\end{equation}
This equation gives a map germ $W: (\mathbb R^2,0) \to \mathbb R^2$ by $W(x,y)= (A(x,y),B(x,y))$.
Now the Chern connection form in the chosen normalization of the forms $\sigma_i$ is
\begin{equation}\label{connection}\textstyle
\gamma=\frac{(2A^2Ax-4A^2By+6ABAy+9BBx)}{4A^3+27B^2}dx+\frac{(4A^2Ay+6ABx+18BBy-9BAx)}{4A^3+27B^2}dy.
\end{equation}
Then the flatness of the web manifest itself as a nonlinear PDE of the second order $d\gamma=0$ (see \cite{Ai} for more detail and the exact form of this PDE).

Consider the variety $E_2$, defined by equation $d\gamma=0$, in the jet space $J^2(\mathbb R^2,\mathbb R^2)$. Total derivatives of this equation with respect to $x$ and $y$  give the {\it prolonged} manifold $E_k$ in each jet space $J^k(\mathbb R^2,\mathbb R^2)$ for any $k\ge 2$.  Differentiating $W$ we obtain a (possibly singular)  parameterized  manifold $\Sigma_k\subset J^k(\mathbb R^2,\mathbb R^2)$ of dimension $\dim \Sigma_k \le 2$.

Conversely, any map $W:(x,y)\to (A(x,y),B(x,y))$ with $\delta:=-4A^3-27B^2\not\equiv 0$ locally determines a 3-web, namely, the web of integral curves of the implicit ODE (\ref{cub}). Any property of the web can be described in terms of the map $W$ and its derivatives.

 In particular, the failure of the following {\it regularity condition} on the criminant $C$
 \begin{equation}\label{regularity}
{\rm rank}((x,y,p)\mapsto (F,F_p))=2,
\end{equation}
 used in \cite{Dn,Ai},
determines a submanifold $R$ of codimension 3 in $J^1(\mathbb R^2,\mathbb R^2)$ by:
\begin{equation}\label{R}
\left\{
\begin{array}{l}
4A^3+27B^2=0\\
A_xB_y-A_yB_x=0\\
3BA_x-2AB_x=0\\
3BA_y-2AB_y=0.
\end{array}
\right.
\end{equation}
One obtains  equations (\ref{R})  by direct computation, substituting the double root $p$ of (\ref{cub}), satisfying $F=F_p=0$, into the determinantal variety determined by the rank drop
$${\rm rank}((x,y,p)\mapsto (F,F_p))<2.$$
Observe that the regularity condition (\ref{regularity}) implies that $M$ and $C$ are  smooth.

Let us call the set of points $m$ on the criminant $C$ a Legendrian locus, if
 the tangent plane to $M$ at $m\in C \subset M$ coincides with the contact plane  ${\rm d} y-p{\rm d}x=0$  at $m$.
The Legendrian locus is determined by  the following conditions:
$$
dF\wedge ({\rm d}y-p{\rm d}x)=0, \ \ \ \ \ \ \ \ \  F=Fp=0.
$$
 Direct computation yields the following
 equations for the corresponding variety in $J^1(\mathbb R^2,\mathbb R^2)$, which we denote by $L$:
\begin{equation}
\left\{
\begin{array}{l}
4A^3+27B^2=0\\
6ABA_x-9B^2A_y+6ABB_y-4A^2B_x=0
\end{array}
\right.
\end{equation}
Note that the variety $L$ is of codimension 2.
\begin{definition}\label{defgen}
 We call a flat web, formed by integral curves of equation (\ref{cub}), generic if the parameterized  manifold $\Sigma_1$ is transverse to the variety $R\cup L$.
\end{definition}

\noindent {\bf Remark.} The varieties $L$ and $R$ are defined in geometric terms. Therefore they are invariant with respect to any diffeomorphisms preserving the form of equation (\ref{cub}). A diffeomorphism preservs this form if it respects  the condition $p_1+p_2+p_3=0$ and does not send any of the roots $p_i$ to infinity.

\section{Generic singularities of  flat 3-webs}
\begin{theorem}\label{classification}
A generic singular flat 3-web germ  defined by implicit cubic ODE (\ref{cub})
 is diffeomorphic either to the web formed by integral curves of the equation
\begin{equation}\label{3}
p^3+2xp+y=0,
\end{equation}
if all 3 web leaves are tangent at the base point, or to the web, formed by the lines $x=cst$ and by the integral curves of the equation
\begin{equation}\label{2}
p^2-x=0,
\end{equation}
if only 2 web leaves are tangent at the base point.
\end{theorem}
{\it Proof:} Let us fix the base point $(x_0,y_0)$ of the web germ and the point $m=(x_0,y_0,p_0)\in C\subset M$. Since the variety $R$ is of codimension 3, the surface $\Sigma_1$ does not intersect $R$. Hence the regularity condition (\ref{regularity}) is satisfied. The variety $L$ is of codimension 2. Therefore the surface $\Sigma_1$ can meet $L$ only in  isolated points.

Consider first the case of a triple root of (\ref{cub}).  As shown in \cite{Ai} (see Lemma 2 there), the contact plane is tangent to $M$ at $m$. (This amounts to $B_x=0$ at $(x_0,y_0)$, since  one has $p_0=0$ for the normalization (\ref{cub}).) Thus $\Sigma_1$ intersects $L$ at a point. Due to genericity, this point is isolated and the criminant $C$ is transverse to the contact plane field in some punctured neighborhood of $m$. Consequently the conditions of Theorem 7 of \cite{Ai} are satisfied and equation (\ref{cub}) is equivalent to the form (\ref{3}) with respect to some local diffeomorphism at $\pi (m)=(x_0,y_0)$.

Now consider the case of double root. Equation (\ref{cub}) have a quadratic factor:
$$
p^3+A(x,y)p+B(x,y)=(p^2+a(x,y)p+b(x,y))(p-a(x,y))=0.
$$
At the point $(x_0,y_0)$ we have a simple root $p_1=a(x_0,y_0)\ne 0$ and a double root  $p_0=-a(x_0,y_0)/2$. Since these roots are distinct, the point $(x_0,y_0,p_0)$ lies on the surface defined by the quadratic equation $f(x,y,p)=p^2+a(x,y)p+b(x,y)=0$. Hence the regularity condition (\ref{regularity}) is now fulfilled for the function $f$.
If $\Sigma_1$ misses the variety $L$, then  the contact plane is not tangent to $C$ at $m$. Due to the regularity
 condition, Theorem 5 of \cite{Ai} implies that our web is equivalent to the normal form (\ref{2}). Note that for this form the discriminant $x=0$ is a leave of one of the web foliations. Moreover,  the discriminant curve $a^2-4b=0$ of $p^2+a(x,y)p+b(x,y)$ is smooth. In fact, the regularity condition implies that at least one of the expressions $b_x+p_0a_x=b_x-\frac{a}{2}a_x$ and $b_y+p_0a_y=b_y-\frac{a}{2}a_y$ does not vanish at  $(x_0,y_0,p_0)$. Thus the gradient of $a^2-4b$ is not zero.

Finally, if  $\Sigma_1$ intersects $L$ at a point, then the  point $m$ of the criminant is Lagrangian.
 As we have mentioned before, the discriminant curve is tangent to the vector field $\partial_x+a\partial_y$ at each point $(x,y)\ne (x_0,y_0)$, therefore this also holds true at $(x_0,y_0)$. But at the point $(x_0,y_0)$  the vector field $\partial_x-\frac{a}{2}\partial_y$ is tangent to the discriminant curve as well, since $m$ is Lagrangian. Thus $a(x_0,y_0)=-\frac{a(x_0,y_0)}{2}$ and the root is triple. This contradiction finishes the proof.
\hfill $\Box$\\

\begin{figure}[th]
\begin{center}
\epsfig{file=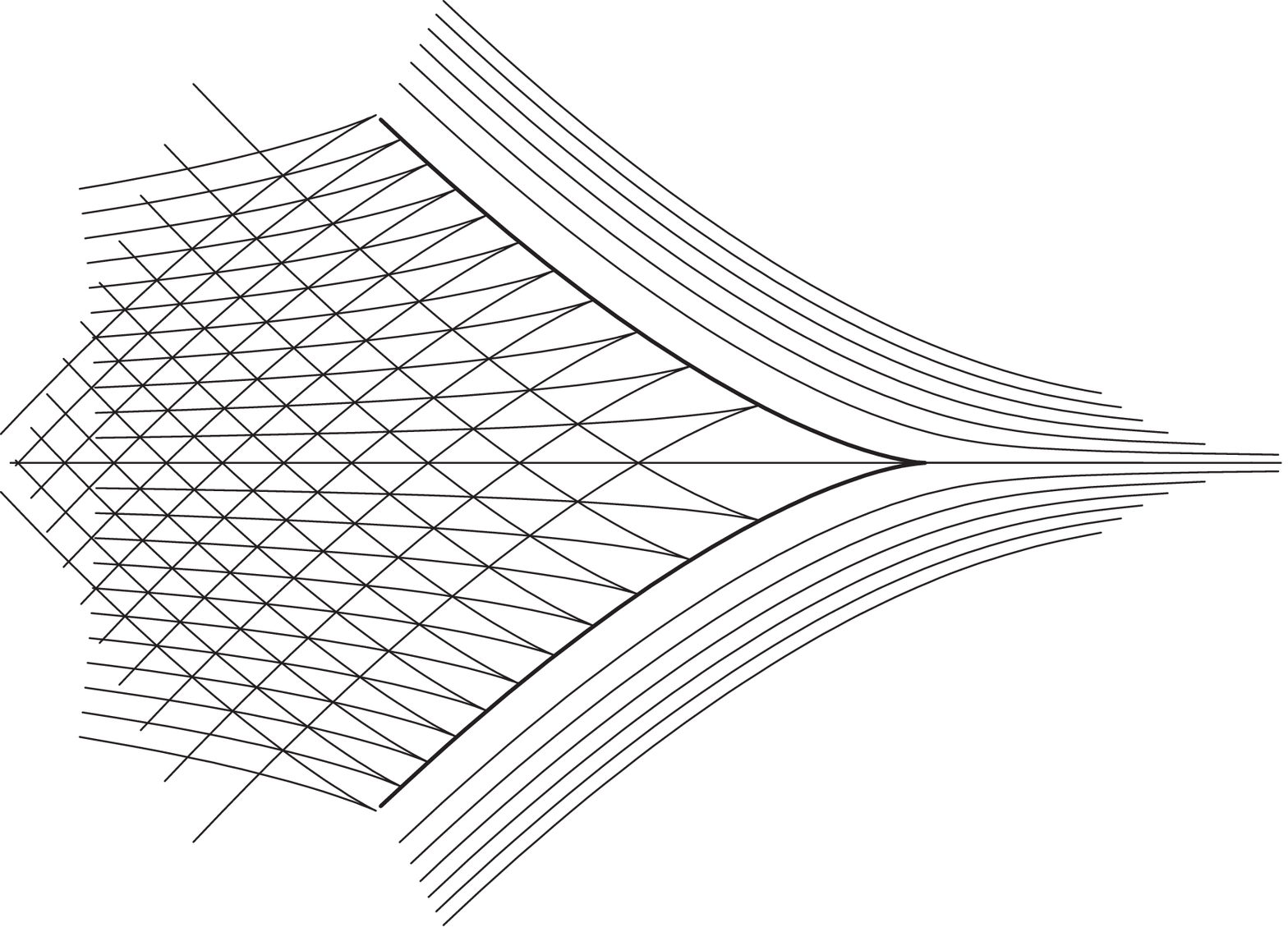,width=40mm, angle=270} \ \ \ \ \ \ \ \ \ \ \ \ \ \ \ \
\ \ \ \ \ \ \
\epsfig{file=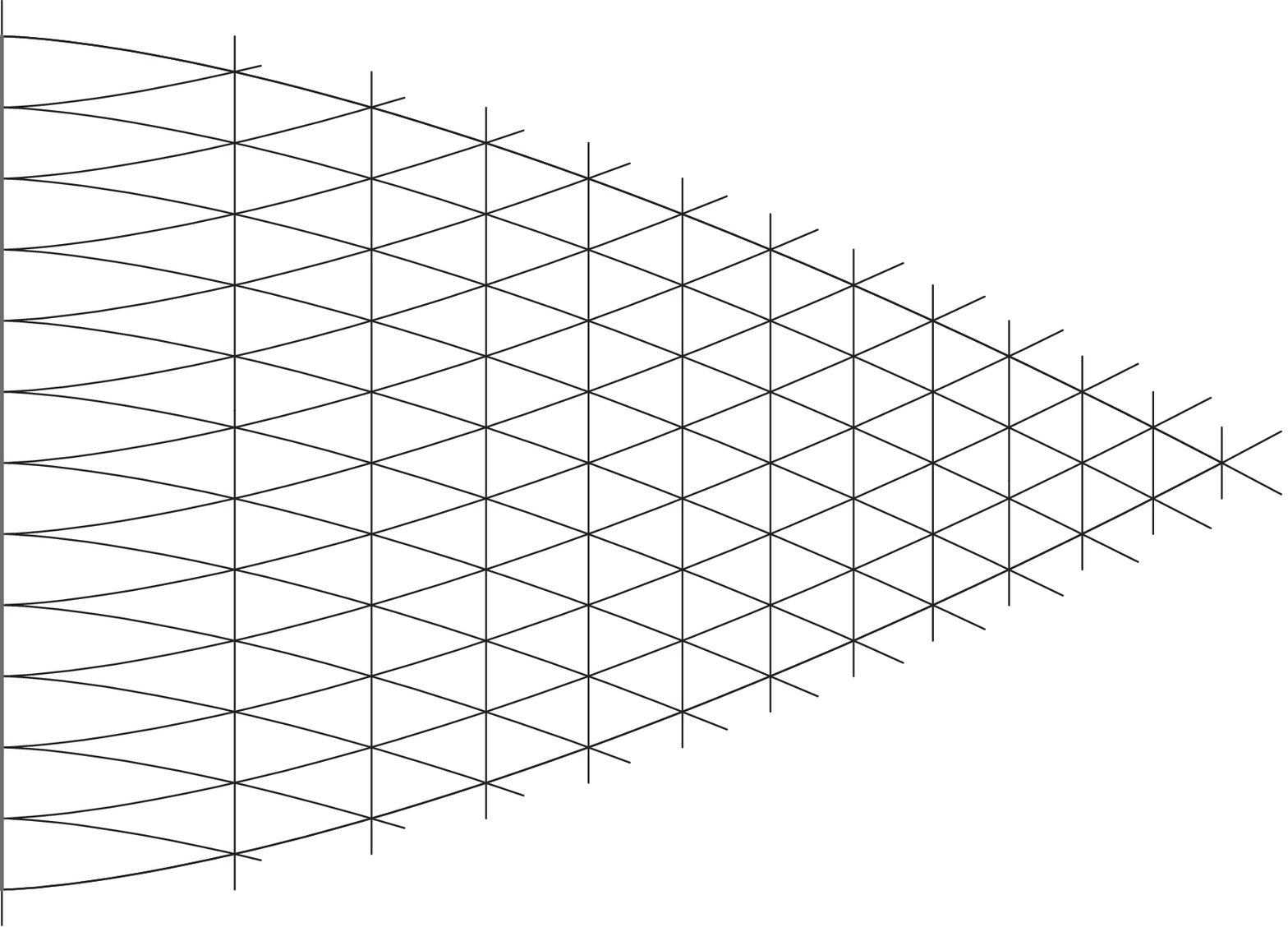,width=35mm, angle=270} \caption{Flat 3-webs
of $p^3+2xp+y=0$ (left) and  of $p^2=x$ (right).} \label{Pic_sol}
\end{center}
\end{figure}

\noindent {\bf Remark 1.} The web germ (\ref{2}) is equivalent to the web germ of the web (\ref{3}) at any point $(x_0,y_0)\ne (0,0)$ on the discriminant. See Figure \ref{Pic_sol}.

\medskip
\noindent {\bf Remark 2.}  As examples of non-generic singular flat 3-webs, one easily constructs the following ones, starting from the well-folded singularity $p^2+ax^2=2y$, $a=cst$:
$$
(xp-2y)(p^2+ax^2-2y)=0, \ \ \ \ \ (2yp-x(2y-ax^2))(p^2+ax^2-2y)=0.
$$
(See  \cite{Dn} for the detailed study of  well-folded singularities of implicit quadratic ODEs.)
The point $(0,0,0)$ is Legendrian for the folded point $(0,0)$ of  $p^2+ax^2=2y$, and all coefficients of the above cubic  equations vanish at the base point.

\section{Concluding remarks}\label{CR}

\noindent $\bullet$ {\bf Genericity and differential consequences of PDE.} It seems natural to expect  a normal form with a triple root, where $\Sigma_1$ misses $L$.
As it was shown in \cite{Ai}, the surface $\Sigma_1$ always meets  $L$ whenever the regularity condition (\ref{regularity}) holds true and all 3 roots $p_i$ of (\ref{cub}) coincide. It is interesting that this differential condition of  order one follows from the PDE $d\gamma=0$, which is of order two.

\noindent $\bullet$ {\bf Perturbations of singular flat 3-webs.}
The variety $L$ is of codimension 2 and the variety $R$ is of codimension 3, thus, due to Thom's Transversality Lemma, for a  generic smooth map $U\subset \mathbb R^2 \to  \mathbb R^2$, its lift to $U\subset \mathbb R^2 \to  J^1(\mathbb R^2,\mathbb R^2)$ misses $R$ and has isolated points of intersection with $L$. But our map $W$ satisfies the PDE $d\gamma =0$ and its differential consequences, which defines an algebraic provariety of {\it infinite codimension} in the jet space $J^{\infty}(\mathbb R^2,\mathbb R^2)$. Therefore $W$ is not generic in the definition of Tougeron \cite{Ti}.

It remains a challenge to show rigorously that for any singular flat 3-web there is a small perturbation, {\it within the class of flat 3-webs}, that  brings the web to a generic flat 3-web in the sense of Definition \ref{defgen}. For the case when $A,B$ are of order 2 the existence of such a perturbation can be shown directly, avoiding the difficult problem of existence of solutions to the singular PDE $d\gamma =0$. In fact, a change of variables $y=F(X,Y),\ x=G(X,Y)$ with $F(0,0)=G(0,0)=0$ preserving the form (\ref{cub}), sends (\ref{3}) to the equation $P^3+\widetilde{A}(X,Y)P+\widetilde{B}(X,Y)=0$ with  $\widetilde{A}_X(0,0)=2G_X^3(0,0)/F_Y^2(0,0)$, $\widetilde{A}_Y(0,0)=5G_X^2(0,0)G_Y(0,0)/F_Y^2(0,0)$, $\widetilde{B}_X(0,0)=0$ and $\widetilde{B}_Y(0,0)=G_X^3(0,0)/F_Y^2(0,0)$. Therefore keeping $F_YG_X\ne 0$, which ensures that $F,G$ are local coordinates, one can make all the first derivatives of $\widetilde{A},\widetilde{B}$ arbitrarily small. Thus, according to our definition of perturbation in Section \ref{section_genericity}, the web of $P^3+\widetilde{A}(X,Y)P+\widetilde{B}(X,Y)=0$ is the desired perturbation.

\noindent $\bullet$ {\bf Flat 3-webs of characteristics of PDEs.} In this note we define a singular web as a web of integral curves of some implicit ODE. This approach is motivated by examples from mathematical physics. For instance, any solution of  the associativity equation
\begin{equation}\label{ass}
u_{yyy}=u_{xxy}^2-u_{xxx}u_{xyy} \end{equation}
determines the so called characteristic 3-web in $(x,y)$-plane  by the following cubic implicit ODE:
$$
\mu ^3+u_{xxx}\mu ^2-2u_{xxy}\mu +u_{xyy}=0,\ {\rm where}  \
\mu=-\frac{dx}{dy}.
$$
This web is flat. Each weighted homogeneous solution of
associativity equation (\ref{ass}) defines a structure of
Frobenius 3-manifold (see \cite{Df}), for example, the web of normal form (\ref{3}) is the characteristic web of a polynomial solutions to (\ref{ass}). (See \cite{Aw,Af} for geometric construction of such webs from Frobenius manifolds.)

Another class of examples of flat 3-webs comes as characteristic webs of integrable systems (see \cite{Fc,Fi})  of hydrodynamic type:
$$
\vec{v}_x=\Lambda(\vec{v})\vec{v}_y,
$$
where  $\vec{v}: U\subset \mathbb R^2 \to \mathbb R^3$ is a vector function and $\Lambda(v)$ is a $3\times 3$ matrix,   whose entries depend explicitly only on $\vec{v}$. Again, each solution $\vec{v}(x,y)$ defines the characteristic web by the following implicit ODE
$$
\det[\Lambda(\vec{v}(x,y))+p\cdot\mathbf{1}]=0.
$$
Here determinantal singularities may come into play and it would be surprising if generic singularities of characteristic webs of scalar PDE and those of systems of hydrodynamic types were the same. (Sometimes a system of hydrodynamic type is equivalent to a scalar PDE, see \cite{Al}.)

\section*{Acknowledgements}

The author thanks J.~Rieger for useful discussions.
This research was partially supported by  DAAD-CAPES grant A/13/03350 and  FAPESP grant \#2014/17812-0.

\end{document}